\documentclass[10pt,a4paper]{amsart}

\usepackage{amssymb}
\usepackage{amsthm}
\usepackage{amsfonts}
\usepackage{microtype}
\usepackage{mathrsfs}
\usepackage{graphicx}
\usepackage{color}
\usepackage{latexsym}
\usepackage{rotating}
\usepackage{xspace}
\usepackage[all]{xy}
\usepackage{longtable}
\usepackage{leftidx}
\usepackage{mathtools}
\usepackage{titlesec}
\usepackage{tikz}
\usetikzlibrary{calc}
\usetikzlibrary{arrows}
\usetikzlibrary{decorations.markings}
\usepackage{geometry}

\newtheorem{theorem}{Theorem}[section]
\numberwithin{equation}{theorem}
\newtheorem{lemma}[theorem]{Lemma}
\newtheorem{proposition}[theorem]{Proposition}
\newtheorem{corollary}[theorem]{Corollary}

\newtheorem{definition}[theorem]{Definition}
\newtheorem{example}[theorem]{Example}

\newtheorem{remark}[theorem]{Remark}
\newtheorem{definition and remark}[theorem]{Definition and Remark}

\newtheorem{definitionandremark}[theorem]{Definition and Remark}
\newtheorem{acknowledgement}{Acknowledgement}

\newcommand{\Ass}{\operatorname{Ass}}

\newcommand{\grade}{\operatorname{grade}}

\newcommand{\Spec}{\operatorname{Spec}}

\newcommand{\ara}{\operatorname{ara}}

\newcommand{\cd}{\operatorname{cd}}
\newcommand{\RH}{\operatorname{H}}

\newcommand{\Ht}{\operatorname{ht}}

\newcommand{\pd}{\operatorname{pd}}

\newcommand{\V}{\operatorname{V}}

\newcommand{\Ext}{\operatorname{Ext}}
\newcommand{\Supp}{\operatorname{Supp}}

\newcommand{\Hom}{\operatorname{Hom}}

\newcommand{\Ann}{\operatorname{Ann}}
\newcommand{\Rad}{\operatorname{Rad}}

\newcommand{\depth}{\operatorname{depth}}

\newcommand{\lo}{\longrightarrow}
\newcommand{\fm}{\frak{m}}
\newcommand{\fp}{\frak{p}}

\newcommand{\fa}{\frak{a}}
\newcommand{\fb}{\frak{b}}

\newcommand{\suchthat}{\;\ifnum\currentgrouptype=16 \middle\fi|\;}

\DeclareMathSymbol{\perp}{\mathrel}{symbols}{"3F}

\newenvironment{prf}[1][Proof]{\begin{proof}[\bf #1]}{\end{proof}}

\makeatletter
\newcommand{\holim@}[2]{%
\vtop{\m@th\ialign{##\cr
\hfil$#1\operator@font holim$\hfil\cr
\noalign{\nointerlineskip\kern1.5\ex@}#2\cr
\noalign{\nointerlineskip\kern-\ex@}\cr}}%
}
\newcommand{\holim}{%
\mathop{\mathpalette\holim@{\rightarrowfill@\textstyle}}\nmlimits@
}
\makeatother

\makeatletter
\def\@secnumfont{\bfseries}
\def\section{\@startsection{section}{1}%
\z@{.7\linespacing\@plus\linespacing}{.5\linespacing}%
{\normalfont\Large\bfseries\filcenter}}
\def\subsection{\@startsection{subsection}{2}%
\z@{.5\linespacing\@plus.7\linespacing}{-.5em}%
{\normalfont\large\bfseries}}
\makeatother
\begin{document}

\author[K. Divaani-Aazar, A. Ghanbari Doust, M. Tousi and Hossein Zakeri]
{Kamran Divaani-Aazar, Akram Ghanbari Doust, Massoud Tousi\\ and\\ Hossein Zakeri}

\title[Modules whose finiteness ...]
{Modules whose finiteness dimensions coincide with their cohomological dimensions}

\address{K. Divaani-Aazar, Department of Mathematics, Faculty of Mathematical Sciences, Alzahra
University, Tehran, Iran-and-School of Mathematics, Institute for Research in Fundamental Sciences (IPM),
P.O. Box 19395-5746, Tehran, Iran.}
\email{kdivaani@ipm.ir}

\address{A. Ghanbari Doust,  Faculty of Mathematical Sciences and Computer, Kharazmi University, Tehran, Iran.}
\email{fahimeghanbary@yahoo.com}

\address{M. Tousi, Department of Mathematics, Faculty of Mathematical Sciences, Shahid Beheshti University, Tehran,
Iran.}
\email{mtousi@ipm.ir}

\address{H. Zakeri, Faculty of Mathematical Sciences and Computer, Kharazmi University, Tehran, Iran.}
\email{hoszakeri@gmail.com}

\subjclass[2010]{13C14; 13C05; 13D45.}

\keywords {Arithmetic rank; cohomological dimension; Koszul complex; local cohomology; relative system of parameters.}

\begin{abstract}  Let $\fa$ be an ideal of a commutative Noetherian ring $R$ with identity. We study finitely generated
$R$-modules $M$ whose $\fa$-finiteness and $\fa$-cohomological dimensions are equal. In particular, we examine relative
analogues of quasi-Buchsbaum, Buchsbaum and surjective Buchsbaum modules. We reveal several interactions between these
types of modules that extend some of the existing results in the classical theory to the relative one.
\end{abstract}

\maketitle

\tableofcontents

\section{Introduction}

Throughout this paper, $R$ is a commutative Noetherian ring with identity, $\fa$ is an ideal of $R$ and $M$ is a finitely
generated $R$-module. The $\fa$-cohomological dimension of $M$, $\cd\left(\fa,M\right)$, is defined as the supremum of
non-negative integers $i$ for which $\text{H}_{\fa}^i(M)\neq 0$. Also, the $\fa$-finiteness dimension of $M$,
$\text{f}_{\fa}\left(M\right)$, is defined as the infimum of non-negative integers $i$ for which $\text{H}_{\fa}^i(M)$ is
not finitely generated. If $\text{f}_{\fa}\left(M\right)<\infty$, then it is known that $\text{f}_{\fa}\left(M\right)\leq
\cd\left(\fa,M\right)$. In this paper, we study the behavior of the $R$-modules $M$ for which the equality holds. The
notion of $\fa$-relative system of parameters plays an essential role in our investigation.

A relative theory of system of parameters is introduced in \cite{DGTZ}. Recall that if $R$ is local with the maximal ideal
$\fm$ and  $\dim M=d$, then a sequence $x_{1}, \ldots, x_{d}$ of elements of $\fm$ forms a system of parameters of $M$ if
and only if the ideals $\langle x_{1},\ldots, x_{d}\rangle$ and $\fm$ have the same radical in $R/\Ann_{R}M$. Although the
Krull dimension seems suitable to define system of parameters, but in the relative case, we appeal to the cohomological dimension
$c:=\cd\left(\fa,M\right)$. Indeed, a sequence $x_{1}, \ldots, x_{c}$ of elements of $\fa$ is said to be an $\fa$-relative
system of parameters, $\fa$-s.o.p, of $M$ if the ideals $\langle x_{1},\ldots ,x_{c} \rangle$ and $\fa$ have the same radical
in $R/\Ann_{R}M$. Contrary to systems of parameters which always exist, relative systems of parameters may don't exist.
It is easy to see that $\fa$ contains an $\fa$-s.o.p of $M$ if and only if $\ara \left(\fa,M\right)=\cd \left(\fa,M\right)$;
see \cite[Lemma 2.2]{DGTZ}. (Here, $\ara \left(\fa,M\right)$ is the arithmetic rank of $\fa$ with respect to $M$ which is
defined as the infimum of the integers $n \in \mathbb{N}_{0}$ for which there exist $x_{1}, \ldots, x_{n}\in R$ such that
the ideals $\langle x_{1}, \ldots, x_{n} \rangle$ and $\fa$ have the same radical in $R/\Ann_{R}M$.)

Suppose that $M\neq 0$. There are many situations in which the equality $\ara \left(\fa,M\right)=\cd \left(\fa,M\right)$
holds. It is clear that if $R$ is local and $\fa$ is the maximal ideal of $R$, then $\ara \left(\fa,M\right)=\cd \left(\fa,
M\right)$. More generally, if $R$ is local and $\fa$ is generated up to radical by a part of a system of parameters of $M$,
then $\ara \left(\fa,M\right)=\cd \left(\fa,M\right)$; see Corollary \ref{4.4}. Also, if $\fa$ is generated by a regular
$M$-sequence, then the equality $\ara \left(\fa,M\right)=\cd \left(\fa,M\right)$ holds. Let $K$ be a field. For a square-free
monomial ideal $\fa$ of a polynomial ring $R=K[x_1,\dots ,x_n]$, it is known that $\cd(\fa,R)=\pd_R\frac{R}{\fa}$; see
\cite[Theorem 1]{Ly}. Characterizing monomial ideals $\fa$ satisfying $\ara(\fa)=\pd_R\frac{R}{\fa}$ has been an active
area of research for years; see e.g. \cite{Ba1}, \cite{Ba2} and \cite{ScV}.

The notion of relative system of parameters enhanced the theory of relative Cohen-Macaulay modules. For instance, when $\fa$
is contained in the Jacobson radical of $R$ and $\ara \left(\fa,M\right)=\cd \left(\fa,M\right)$,  it is shown that $M$ is
$\fa$-relative Cohen-Macaulay if and only if every $\fa$-s.o.p of $M$ is an $M$-regular sequence if and only if there exists
an $\fa$-s.o.p of $M$ which is an $M$-regular sequence; see \cite[Theorem 3.3]{DGTZ}.

There are several extensions of the notion of Cohen-Macaulay modules. Perhaps, the best known, and geometrically most important,
is that of generalized Cohen-Macaulay modules which is well studied in \cite{T2}. Surjective Buchsbaum, Buchsbaum and quasi
Buchsbaum modules are special types of generalized Cohen-Macaulay modules that have been studied extensively in late 80's; see
e.g. \cite{StV}. The notions of $\fa$-filter regular sequences, $\fa$-weak sequences and $d$-sequences play an important role
in commutative algebra. When $(R,\fm)$ is local, Buchsbaum, quasi Buchsbaum or generalized Cohen-Macaulay $R$-modules can be
characterized by the property that every (some of) their systems of parameters being a $\fm$-weak sequence or $\fm$-filter regular
sequence on them. We intend to characterize the $R$-modules $M$ that every (some of) $\fa$-s.o.p of $M$ is a
$\fa$-weak sequence or $\fa$-filter regular sequence on $M$. This leads to a theory of $\fa$-relative Buchsbaum, $\fa$-relative
quasi Buchsbaum and $\fa$-relative generalized Cohen-Macaulay modules; see Definitions \ref{3.11} and \ref{4.11}.

Next, we will describe the organization of the paper together with its main results.

Section $2$ is devoted to some preliminaries on properties of filter regular sequences, weak sequences and $d$-sequences. In Remark
\ref{2.3}, we reveal some of the existing relationships in the literature between these three types of sequences. Then, we present
some basic characterizations of these sequences. Theorems \ref{2.6} and \ref{2.7} are the main results of this section. The Koszul
complexes, which assist us through the paper, are the main ingredients of this section.

In Section $3$, we introduce the notion of $\fa$-relative generalized Cohen-Macaulay modules, and establish some of the main
properties of these modules. Cuong, Schenzel and Trung in \cite[Satz 3.3]{CST} showed that every system of parameters of a
generalized Cohen-Macaulay $R$-module $M$ is a filter regular sequence on $M$. Now, the question arises whether every
$\fa$-s.o.p of an $\fa$-relative generalized Cohen-Macaulay $R$-module $M$ is an $\fa$-filter regular sequence on $M$?
Among other things, Theorem \ref{3.3} gives a positive answer to this question.

In Section $4$, we introduce and study $\fa$-relative surjective Buchsbaum, $\fa$-relative Buchsbaum and $\fa$-relative quasi
Buchsbaum modules. Assume that the ideal $\fa$ is contained in the Jacobson radical of $R$. In Theorem \ref{4.1}, we provide a
characterization of $\fa$-relative quasi Buchsbaum modules via local cohomology modules. Theorem \ref{4.2}, among other things,
indicates that a given finitely generated $R$-module is $\fa$-relative Buchsbaum if there exists a generating set of $\fa$
satisfying some certain properties. Theorem \ref{4.8} reveals a close connection between $\fa$-relative Buchsbaum and
$\fa$-relative quasi Buchsbaum $R$-modules.

\section{Koszul complexes and filter regular sequences}

The Koszul complexes play an important role in this paper. When discussing Koszul complexes, we obey the notation of \cite[5.2]{BS},
which shall be used without further comments.

Let $i\in \mathbb{N}_0$. Let $\underline{x}:=x_1,..., x_n$ be a sequence of elements of $R$ and $M$ an $R$-module. It is convenient
for us to set $K^{i}\left(\underline{x},M\right):=K_{n-i}\left(\underline{x},M\right)$. Then the induced co-complex is denoted by $K^{\bullet}\left(\underline{x}, M\right)$, and one has $$\RH^{i}\left(\underline{x},M\right):=\RH^{i}\left(K^{\bullet}
\left(\underline{x}, M\right)\right)=\RH_{n-i}\left(\underline{x}, M\right).$$

As in \cite[Lemma 5.2.2]{BS}, for natural numbers $u,v$ with $u\leq v$, one can define an explicit $R$-homomorphism $$\left(\psi_{u}^{v}\right)_{p}: K_{p}\left(x_{1}^{u},\ldots, x_{n}^{u},M\right)\longrightarrow K_{p}\left(x_{1}^{v},
\ldots, x_{n}^{v},M\right)$$ such that the chain maps $$\left(\psi_{u}^{v}\right)_{\bullet}: K_{\bullet}\left(x_{1}^{u},
\ldots, x_{n}^{u},M\right)\longrightarrow K_{\bullet}\left(x_{1}^{v},\ldots, x_{n}^{v},M\right)$$ turn the family $\{K_{\bullet}\left(x_{1}^{u},\ldots, x_{n}^{u},M\right)\}_{u\in \mathbb{N}}$ into a direct system of complexes and
chain maps. Moreover, it is proved in \cite[Theorem 5.2.9]{BS} that there is a natural $R$-isomorphism $$\RH_{\fa}^{i}
\left(M\right)\cong \lim _{\overrightarrow{u \in \mathbb{N}}}\RH_{n-i}\left(x_{1}^{u}, \ldots, x_{n}^{u} ,M\right).$$
So, there is a canonical $R$-homomorphism $\lambda_{M}^{i}: \RH^{i}\left(\underline{x},M\right)\longrightarrow
\RH_{\fa}^{i}\left(M\right)$.

\begin{lemma}\label{2.1} (See \cite[Lemma 1.5]{StV}.) Let $\fa=\langle x_{1},\ldots, x_{n} \rangle$ and $\fb=\langle y_{1},\ldots,
y_{m} \rangle$ be two ideals of $R$ such that $\fa \subseteq \fb$.  Then for every $R$-module $M$ and every non-negative integer $i$,
there is a commutative diagram
\begin{displaymath}
\xymatrix{\RH^{i}\left(y_{1},\ldots, y_{m}, M\right) \ar[r]  \ar[d] &
\RH^{i}\left(x_{1},\ldots, x_{n},M\right) \ar[d] \\
\RH_{\fb}^{i}\left(M\right) \ar[r] & \RH_{\fa}^{i}\left(M\right).}
\end{displaymath}
\end{lemma}

Next, we recall the definitions of some sequences which play key roles in the study of the modules that will be introduced in
Definitions \ref{3.11} and \ref{4.11}.

\begin{definition}\label{2.2} Let $\fa$ be an ideal of $R$ and  $M$ a finitely generated $R$-module.
\begin{enumerate}
\item[(i)]  Following \cite{CST}, a sequence $x_{1},\ldots, x_{r}$ of elements of $R$ is called an {\it $\fa$-filter regular}
sequence on $M$ if  $$ \langle x_{1},\ldots, x_{i-1} \rangle M:_{M}x_{i} \subseteq \bigcup\limits_{t\in \mathbb{N}}\langle
x_{1},\ldots, x_{i-1}\rangle M:_{M} \fa ^{t}$$ for all $i=1,\ldots, r$.  An $\fa$-filter regular sequence $x_{1},\ldots, x_{r}$
on $M$ is called an \it{unconditioned $\fa$-filter regular} sequence on $M$, if it is an $\fa$-filter regular sequence on $M$
in any order.

\item[(ii)] Following \cite[page 39]{T1}, a sequence $ x_{1},\ldots, x_{r} $ of elements of $R$ is called an  {\it $\fa$-weak
sequence} on $M$ if $$ \langle  x_{1},\ldots, x_{i-1} \rangle M:_{M}x_{i} \subseteq  \langle x_{1}, \ldots, x_{i-1}
\rangle M:_{M}\fa$$ for all $i=1, 2, \ldots, r$. It is an \it{unconditioned $\fa$-weak sequence} on $M$ if $x_{1}^{\alpha_{1}},
\ldots ,x_{r}^{\alpha_{r}}$ is an $\fa$-weak sequence on $M$ in any order for all positive integers $ \alpha_{1},\ldots,
\alpha_{r}$.

\item[(iii)] Following  \cite{T1}, a sequence $x_{1}, \ldots, x_{r}$ of elements of $R$ is called a {\it $d$-sequence} on $M$ if
$$\langle x_{1}, x_2, \ldots, x_{i-1}\rangle M:_{M}x_{i}x_{j}=\langle x_{1}, x_2, \ldots, x_{i-1}\rangle M:_{M}x_{j}$$ for any
$1\leq i \leq j \leq r$. When $x_{1}^{n_{1}}, x_{2}^{n_{2}}, \ldots, x_{r}^{n_{r}}$ form a $d$-sequence on $M$ for all integers
$n_{1}, n_2, \ldots, n_{r}\in \mathbb{N}$, then $x_{1},\ldots, x_{r}$ is called a \it{strong $d$-sequence}  on $M$. A  $d$-sequence on
$M$ is termed \it{unconditioned}, $u.s.d$-sequence, when it forms a strong $d$-sequence on $M$ in any order.
\end{enumerate}
\end{definition}

The above Trung's definition of $d$-sequences is actually a slight weakening of Huneke's definition in \cite{H}. It removes the
assumption that each $x_i$ is not in the ideal generated by the rest of the $x_j$.

In the following remark, we collect some properties of the sequences that are already defined above.

\begin{remark}\label{2.3} Let $\fa$ be an ideal of $R$ and  $M$ a finitely generated $R$-module.
\begin{enumerate}
\item[(A)]
\begin{enumerate}
\item[(i)]  Clearly, every $\fa$-weak sequence on $M$ is an $\fa$-filter regular sequence on $M$. By \cite[Theorem 1.1(vi)]{T1}
a $d$-sequence $x_{1},\ldots ,x_{r}$ on $M$ is an $\fa$-weak sequence on $M$ if $\fa\subseteq \langle x_{1},
\ldots, x_{r}\rangle .$
\item[(ii)] By \cite [Theorem 1.1(vii)]{T1}, a $d$-sequence $x_{1},\ldots, x_{r}$ on $M$ is an  $\fa$-filter regular sequence
on $M$ if $\fa \subseteq \Rad\left(\langle x_{1},\ldots, x_{r} \rangle \right).$ Also by \cite[Proposition 2.1]{T1}, if
$x_{1},\ldots, x_{r} \in \fa$ is an $\fa$-filter regular sequence on $M$, then there exist natural numbers $ s_{1}\leq \cdots
\leq s_{r} $ such that $x_{1}^{s_{1}}, \ldots, x_{r}^{s_{r}}$ is a $d$-sequence on $M$.
\item[(iii)] By \cite[Proposition 2.2]{T1}, $ x_{1},\ldots ,x_{r} \in \fa  $ is a $ d $-sequence on $M$ if one of the following
conditions is satisfied:
\begin{enumerate}
\item[(a)] $x_{1},\ldots, x_{i-1},x_{i}^{2}$ is an $\fa$-weak sequence on $M$ for all $ i=1,\ldots, r$.
\item[(b)] $x_{1},\ldots, x_{r}$ is an $\fa$-weak sequence on $M$ in $\fa ^{2}$.
\end{enumerate}
\end{enumerate}
\end{enumerate}
\begin{enumerate}
\item[(B)] Below, we record three essential properties of filter regular sequences.
\begin{enumerate}
\item[(i)] Let  $\underline{x}=x_{1},\ldots, x_{r} \in \fa$  be an $\fa$-filter regular sequence on $M$. Then by \cite[Proposition 2.3]{AS},
there are the natural isomorphisms
\begin{equation*}
\RH_{\fa}^{i}\left(M\right) \cong
\begin{cases}
\RH_{\langle \underline{x} \rangle }^{i}\left(M\right) &  \text{for}\quad 0 \leq i <r\\
\RH_{\fa}^{i-r}\left(\RH_{\langle \underline{x} \rangle}^{r}\left(M\right)\right) & \text{for} \quad i \geq r.
\end{cases}
\end{equation*}
\item[(ii)] Let $\fa=\langle x_{1}, \ldots, x_{r} \rangle$. By \cite[Proposition 1.2]{TZ}, $\fa$ has generators $ y_{1},\ldots,
y_{r}$ which form an unconditioned $\fa$-filter regular sequence on $M$. We call a such set of generators for $\fa$, an
$\text{f}$-generating set of $\fa$ with respect to $M$. Note that for any positive integer $\ell$, we may find an $\text{f}$-generating
set of $\fa$ with respect to $M$ such that it has more than $\ell$ elements.
\item[(iii)]  A sequence $x_{1},\ldots, x_{r}$ of elements of $R$ is an $\fa$-filter regular sequence on $M$ if  and only if for every
$1\leq i\leq r$, the element $x_{i}$ doesn't belong to the union of the elements of $\Ass_{R}\left(\frac{M}{\langle x_{1}, \ldots, x_{i-1} \rangle M}\right) \setminus
\V\left(\fa\right)$.
\end{enumerate}
\end{enumerate}
\end{remark}

In the next result, we consider the situation in which a given filter regular sequence on $M$ forms a weak sequence on $M$.

\begin{lemma}\label{2.4} Let $\fa$ be an ideal of $R$ and $M$ a finitely generated $R$-module. Let $\underline{x}:=x_{1},\ldots,
x_{n}$ be a sequence of elements of $\fa$ which forms an unconditioned $\fa$-filter regular sequence on $M$. Assume that the
natural map $\lambda_M^{i}:\RH^{i}\left(\underline{x},M\right)\longrightarrow \RH_{\langle \underline{x} \rangle}^{i}\left(M\right)$
is surjective for all $0 \leq i \leq n-1$. Then $\underline{x}$ is an unconditioned $\langle \underline{x}\rangle $-weak sequence on
$M$. In particular, $\underline{x}$ is an  $u.s.d$-sequence on $M$.
\end{lemma}

\begin{prf} Set $\fb:=\langle \underline{x} \rangle$. By induction on $0 <i \leq n $, we show that for any permutation $\delta$ of the set
$\{1,2, ..., n\}$ and all positive integers $t_1, \ldots, t_n$, the sequence $x_{\delta (1)}^{t_1},\ldots, x_{\delta (i)}^{t_i}$ is a $\fb$-weak
sequence on $M$. Since $\lambda_M^{0}$ is surjective, $\left(0:_{M} \fb\right)=\RH_{\fb}^{0}\left(M\right)$. By Remark \ref{2.3}(B)(i), one has:
$$\RH_{\langle x_{\delta (1)}^{t_1} \rangle }^{0}\left(M\right)=\RH_{\fa }^{0}\left(M\right)=\RH_{\fb }^{0}\left(M\right),$$  and so
$$\left(0:_{M}x_{\delta (1)}^{t_1}\right) \subseteq \RH_{\langle x_{\delta (1)}^{t_1} \rangle }^{0}\left(M\right)=\RH_{\fb }^{0}\left(M\right)=
\left(0:_{M}\fb\right). $$ Hence, the case $i=1$ holds.

Now, let $2\leq i \leq n$ and the result has been proved for $i-1$. Let $\sigma$ be a permutation of the set $\{1,2, ..., n\}$ and $u_1, \ldots,
u_n$ be positive integers. Set $z_i:=x_{\sigma (i)}^{u_i}$ for every $1\leq i\leq n$. By the induction hypothesis, $z_1,\ldots, z_{i-1}$ is a $\fb$-weak
sequence on $M$, and it remains to show that $$\langle  z_1,\ldots, z_{i-1} \rangle M:_{M}z_{i} \subseteq  \langle z_1,\ldots, z_{i-1} \rangle M:_{M}\fb.$$
By \cite[Lemma 3.5]{T2}, one deduces that $$\langle z_{1}^{k},\ldots, z_{i-1}^{k}\rangle M:_{M} \left(z_{1}\ldots  z_{i-1}\right)^{k-1}=\langle
z_{1},\ldots, z_{i-1}\rangle M+\sum_{j=1}^{i-1} \left(\langle z_{1},\ldots ,\hat{z}_{j} ,\ldots, z_{i-1}\rangle  M:_{M} \fb \right)$$ for all $k\geq 2$.
On the other hand, we have the following obvious containment $$\langle z_{1},\ldots, z_{i-1}\rangle M+\sum_{j=1}^{i-1} \left(\langle z_{1},\ldots,
\hat{z}_{j}, \ldots, z_{i-1}\rangle  M:_{M} \fb\right) \subseteq \langle z_{1},\ldots, z_{i-1} \rangle M:_{M} \fb. $$ So the argument will be complete,
if we show that for each $m_{i}\in \langle z_{1},\ldots, z_{i-1} \rangle M:_{M}z_{i}$, there exist $n_{i}\in \langle z_{1},\ldots, z_{i-1} \rangle M:_{M}
\fb$ and $l\geq 2$ such that $$m_{i}-n_{i} \in \langle z_{1}^{l},\ldots, z_{i-1}^{l} \rangle M:_{M} \left(z_{1}\ldots z_{i-1}\right)^{l-1}.$$

Let $m_{i}\in \langle z_{1},\ldots, z_{i-1} \rangle M:_{M}z_{i}$. Then $z_{i}m_{i}=-\sum_{j=1}^{i-1}z_{j}m_{j}$ for some elements $m_1,\ldots, m_{i-1}$
in $M$. Clearly $t:=\sum_{j=1}^{i}m_{j}e_{j} \in \text{ker}d_{1}$, where $d_{1}:K_{1}\left(z_{1},\ldots, z_{i},M\right) \longrightarrow K_{0}
\left(z_{1},\ldots, z_{i},M \right)$ is the first differential map of the Koszul complex $K_{\bullet}\left(z_{1},\ldots, z_{i},M\right)$.
By Lemma \ref{2.1}, the following diagram is commutative:
\begin{displaymath}
\xymatrix{\RH^{i-1}\left(\underline{x},M\right) \ar[r]^{\theta}  \ar[d]_{\lambda_M^{i-1}} & \RH^{i-1}\left(z_{1},\ldots, z_{i}, M\right)
\ar[d]_{g_{1}} \\
\RH_{\fb}^{i-1}\left(M\right) \ar[r]^{h} & \RH_{\langle z_{1},\ldots, z_{i} \rangle}^{i-1}\left(M\right)}
\end{displaymath}
By Remark \ref{2.3}(B)(i), $h$ is an isomorphism. Also by the assumption, $\lambda_M^{i-1}$ is surjective. So, there is $z \in \RH^{i-1}
\left(\underline{x},M\right)$ such that $g_{1}\left(\theta \left(z\right)\right)=g_{1}\left(t+\text{im}d_{2}\right)$. Suppose that $\theta
\left(z\right)=f+\text{im}d_{2}$, where $f=\sum_{j=1}^{i}n_{j}e_{j} $. Since $\fb \RH^{i-1}\left(\underline{x},M\right)=0 $, it follows
that $\fb f \subseteq \text{im}d_{2}$. It implies that $x_{j}n_{i} \in \langle z_{1}, \ldots, z_{i-1} \rangle M $ for all $j=1,\ldots, n$,
and so $n_{i} \in \langle z_{1},\ldots, z_{i-1}\rangle M :_{M} \fb$. Now, one has $$ 0=g_{1}\left(\left(t-f\right)+\text{im}d_{2}\right)\in
\RH_{\langle  z_{1},\ldots, z_{i} \rangle}^{i-1}\left(M\right)\cong \lim_{\overrightarrow{\alpha}}\RH_{1}\left(z_{1}^{\alpha},\ldots,
z_{i}^{\alpha},M\right).$$ So, there exists an integer $l\geq 2$ such that $\overline{\left(\psi^{l}_{1}\right)_{1}} \left(\left(t-f\right)+\text{im}d_{2}\right)=0$, where $$\overline{\left(\psi_{1}^{l}\right)_{1}}:\RH_{1}\left(z_{1},\ldots, z_{i},M\right) \longrightarrow \RH_{1}\left(z_{1}^{l},\ldots, z_{i}^{l},M\right)$$ is the natural map induced by the map $\left(\psi_{1}^{l}\right)_{1}$
given in the paragraph preceding Lemma \ref{2.1}. Thus, we have $$\left(z_{1}\ldots z_{i-1}\right)^{l-1}\left(m_{i}-n_{i}\right)\in \langle z_{1}^{l},\ldots, z_{i-1}^{l}\rangle M.$$

For the second assertion, note that by Remark \ref{2.3}(A)(iii), every unconditioned $\langle \underline{x} \rangle$-weak sequence on $M$ is
an $u.s.d$-sequence on $M$. So, it is obvious by the first assertion.
\end{prf}

The following lemma together with Theorem \ref{2.6} is required to prove Theorem \ref{2.7}.

\begin{lemma}\label{2.5} Let $\underline{x}:=x_{1},\ldots, x_{s}$ be a sequence of elements of $R$ and set $\fa:=\langle \underline{x} \rangle$.
Let $M$ be a finitely generated $R$-module and $n\leq s$ a natural number. Suppose that every $n$ elements of $\lbrace x_{1}, \ldots, x_{s}
\rbrace$ forms a $d$-sequence on $M$. Then the natural  $R$-homomorphism $\RH_{i}\left(\underline{x},0:_{M}\fa\right)\longrightarrow \RH_{i}\left(
\underline{x},M\right)$ induced by the inclusion map $0:_{M}\fa \hookrightarrow M$ is injective for every nonnegative integer $i$
that is satisfying $0 \leq s-i \leq n$.
\end{lemma}

\begin{prf}  Let $i$ be a nonnegative integer such that $0 \leq s-i \leq n$ and set $K_{i}:=K_{i}\left(\underline{x},M\right) $ and
$\overline{K_{i}}:=K_{i}\left(
\underline{x},0:_{M}\fa\right)$. Let $\RH_{i}$ and $\overline{\RH_{i}}$ denote the $i$th homology modules of these complexes;
respectively. Since
$x_{j}\left(0:_{M}\fa\right)=0 $ for all $j=1, \ldots, s$, one has $\overline{\RH_{i}}=\overline{K_{i}}$. We have to show that the
natural $R$-homomorphism
$$\overline{K_{i}} \longrightarrow \RH_{i}\left(\underline{x},M\right)$$ $$z\mapsto z+\text{im}d_{i+1}$$ is injective.
That is $\overline{K_{i}} \bigcap \text{im}d_{i+1}=0.$  Let $x \in \overline{K_{i}} \bigcap \text{im}d_{i+1}$. Then $$ x=\sum_{ l \in
\mathcal{I}\left(i,s\right)}
\acute{m}_{l_{1}\ldots l_{i}}e_{l_{1}} \ldots e_{l_{i}}  \in \overline{K_{i}} $$ and there is $$ y=\sum_{ j \in \mathcal{I}\left(i+1 ,
s\right) } m_{j_{1} j_{2}
\ldots j_{i+1}}e_{j_{1}}e_{j_{2}}\ldots e_{j_{i+1}} \in K_{i+1} $$ such that $d_{i+1}\left(y\right)=x$. So,
$$\sum_{ l \in \mathcal{I}\left(i,s\right) }\acute{m}_{l_{1}l_{2} \ldots l_{i}}e_{l_{1}}e_{l_{2}}\ldots e_{l_{i}}=\sum_{j \in \mathcal{I}
\left(i+1, s\right)} \left(\sum _{p=1}^{i+1}\left(-1\right)^{p-1}x_{j_{p}}m_{j_{1}j_{2}\ldots j_{i+1}}e_{j_{1}}\ldots \hat{e}_{j_{p}}
\ldots e_{j_{i+1}}\right).$$ Assume that $k \in \mathcal{I}\left(i,s\right) $ and set $$B_k:=\lbrace \ell \in \lbrace 1,\ldots, s \rbrace|\
k  \  \text{can be induced by deleting}\  \ell   \ \text{from an element of} \ \mathcal{I}\left(i+1 ,s\right) \rbrace.$$ So, $\vert B_k \vert
=s-i$. Hence $$\acute{m}_{k_{1}k_{2}\ldots k_{i}} \in x_{j _{1}}M+x_{j _{2}}M+\ldots +x_{j_{s-i}}M,$$ where $j_{1}, j_{2}, \ldots, j_{s-i}
\in B_k.$ As $x_{j_{1}}\acute{m}_{k_{1}k_{2}\ldots k_{i}}=0$,  \cite[Theorem 1.1]{T1} yields that $$\acute{m}_{k_{1}k_{2}\ldots k_{i}} \in
\left(0:_{M}x_{j_{1}}\right)\bigcap \langle x_{j_{1}}, x_{j_{2}},\ldots, x_{j_{s-i}}\rangle M=0.$$ (Note that, by \cite[Lemma 2.2]{G}, if $a_{1},
\ldots, a_{n}$ is a $d$-sequence on $M$, then, for each $1\leq t \leq n $, $a_{1},a_{2}, \ldots, a_{t}$ is also a $d$-sequence on $M$.)
\end{prf}

The next result connects the notion of weak sequences to that of local cohomology.

\begin{theorem}\label{2.6}  Let $\fa$ be an ideal of $R$ and $M$ a finitely generated $R$-module. Let $x_{1},\ldots, x_{n}\in \fa^{2}$.
Then the following are equivalent:
\begin{enumerate}
\item[(i)] $x_{1},\ldots, x_{n}$ is an $\fa$-weak sequence on $M$.
\item[(ii)] $\fa \RH_{\fa}^{i}\left(M\right)=0$ for all $0 \leq i <n$ and $x_{1},\ldots, x_{n}$ is an $\fa$-filter regular sequence on
$M$.
\end{enumerate}
\end{theorem}

\begin{prf} (i)$ \Longrightarrow $(ii)  Suppose that $x_{1},\ldots, x_{n}$ is an $\fa$-weak sequence on $M$. So, $x_{1},\ldots, x_{n}$
is an $\fa$-filter regular sequence on $M$. Note that by Remark \ref{2.3}(A)(iii), $x_{1},\ldots, x_{n}$ is a $d$-sequence on $M$. Hence,
by Remark \ref{2.3}(B)(i), $$\Gamma_{\fa}\left(M\right)=\Gamma_{ \langle x_{1} \rangle }\left(M\right)=0:_{M}x_{1}=0:_{M}\fa.$$ We use
induction on $n$.  For $n=1$,  we are done. Let $n>1$ and the result has been proved for $n-1$. As $x_{1},\ldots, x_{n-1}$ is an
$\fa$-weak sequence on $M$, the induction hypothesis implies that $\fa\RH_{\fa}^{i}\left(M\right)=0$ for all $0\leq i<n-1$. Thus, it
remains to show that $\fa \RH_{\fa}^{n-1}\left(M\right)=0$. We have the following two exact sequences $$0\longrightarrow 0:_{M}x_{1}\hookrightarrow M \stackrel{\rho}\longrightarrow x_{1}M \longrightarrow 0$$ and $$0 \longrightarrow x_{1}M \stackrel{\lambda}\longrightarrow M \twoheadrightarrow M/x_{1}M \longrightarrow 0, $$
in which all maps are natural. Since $0:_{M}x_{1}=\Gamma_{\fa}\left(M\right) $, we get $\RH_{\fa}^{i}\left(0:_{M}x_{1}\right)=0$ for all $i\geq 1$.
Hence, $\RH_{\fa}^{n-1}\left(\rho\right)$ is an isomorphism. It induces the given isomorphism in the following display:
\begin{align*}
\text{ker}\left(\RH_{\fa}^{n-1}\left(\lambda\right)\right)&\cong \text{ker}\left(\RH_{\fa}^{n-1}\left(\lambda\right) \RH_{\fa}^{n-1}\left(\rho\right)\right)\\
&=\text{ker}\left(\RH_{\fa}^{n-1}\left(\lambda  \rho \right)\right)\\
&=\text{ker}\left(\RH_{\fa}^{n-1}\left(x_{1}1_{M}\right)\right)\\
&=0:_{\RH_{\fa}^{n-1}\left(M\right)}x_{1}.
\end{align*}
As $x_{2},\ldots, x_{n}$ is an $\fa$-weak sequence on $M/x_{1}M$, the induction hypothesis yields that $$\fa \RH_{\fa}^{n-2}\left(M/x_{1}M\right)=0.$$
The exact sequence $$\RH_{\fa}^{n-2}\left(M/x_{1}M\right) \longrightarrow \RH_{\fa}^{n-1}\left(x_{1}M\right) \stackrel{\RH_{\fa}^{n-1}\left(\lambda\right)}\longrightarrow \RH_{\fa}^{n-1}\left(M\right)$$ implies the exact
sequence $$\RH_{\fa}^{n-2}\left(M/x_{1}M\right)\longrightarrow \text{ker}\left(\RH_{\fa}^{n-1}\left(\lambda\right)\right)
\longrightarrow 0,$$ and consequently $\fa\  \text{ker}\left(\RH_{\fa}^{n-1}\left(\lambda\right)\right)=0$.
So, $\fa\left(0:_{\RH_{\fa}^{n-1}\left(M\right)}x_{1}\right)=0 $. Now, one has $$0:_{\RH_{\fa}^{n-1}\left(M\right) }x_{1}
\subseteq 0:_{\RH_{\fa}^{n-1}\left(M\right) } \fa \subseteq 0:_{\RH_{\fa}^{n-1}\left(M\right) } \fa^2 \subseteq
0:_{\RH_{\fa}^{n-1}\left(M\right) }x_{1}.$$ Hence, $0:_{\RH_{\fa}^{n-1}\left(M\right)}\fa=0:_{\RH_{\fa}^{n-1}\left(M\right)}\fa^2$.
It implies that $$\RH_{\fa}^{n-1}\left(M\right)=\bigcup_{t\geq 1} 0:_{\RH_{\fa}^{n-1}\left(M\right) } \fa^t=0:_{\RH_{\fa}^{n-1}
\left(M\right) }\fa,$$ and so $\fa \RH_{\fa}^{n-1}\left(M\right) =0$.

(ii)$\Longrightarrow $(i) We use induction on $n$. Let $n=1$. As $x_{1}$ is an $\fa$-filter regular sequence on $M$,
Remark \ref{2.3}(B)(i) yields that $$0:_{M}x_{1}\subseteq \Gamma_{\langle x_{1} \rangle}\left(M\right)=\Gamma_{\fa}
\left(M\right).$$ Hence, $$\fa\left(0:_{M}x_{1}\right)\subseteq \fa \ \Gamma_{\fa}\left(M\right)=0.$$
Thus,  $x_{1}$ is an $\fa$-weak sequence on $M$.

Next,  let $n>1$ and the result has been proved for $ n-1 $. As $x_{1}$ is an $\fa$-filter regular sequence on $M$ in
$\fa^2$,  \cite[Theorem 1.1]{CQ} implies that  $$\RH_{\fa}^{i}\left(M/x_{1}M\right) \cong \RH_{\fa}^{i}\left(M\right)\bigoplus \RH_{\fa}^{i+1}\left(M\right)$$ for all  $i<n-1$. So, our assumption implies that $\fa \RH_{\fa}^{i}\left(M/x_{1}M\right)=0$
for all  $i<n-1$. Thus by the induction hypothesis, $x_{2}, \ldots, x_{n}$ is an $\fa$-weak sequence on $M/x_{1}M$. As
$x_{1}$  is an $\fa$-weak sequence on $M$, we deduce that $x_{1}, \ldots, x_{n}$ is an $\fa$-weak sequence on $M$.
\end{prf}

\begin{theorem}\label{2.7} Let $\fa$ be an ideal of $R$, $M$ a finitely generated $R$-module and $n$ a natural number. Then the following are equivalent:
\begin{enumerate}
\item[(i)]  Any $n$ elements of every $\text{f}$-generating set $\{x_1, \ldots, x_s\}$ of $\fa$ with respect to $M$ with $s\geq n$ form an unconditioned
$\fa$-weak sequence on $M$.
\item[(ii)] There exists an $\text{f}$-generating set $\{x_1, \ldots, x_s\}$ of $\fa$ with respect to $M$ such that $s\geq n$ and any $n$ elements in
it form an unconditioned $\fa$-weak sequence on $M$.
\item[(iii)] There exists an $\text{f}$-generating set $\{x_{1},\ldots, x_{s}\}$ of $\fa$ with respect to $M$ with $s\geq n$ such that the natural map
$\lambda_{M}^{i}:\RH^{i}\left(x_{1}, \ldots, x_{s},M\right) \longrightarrow \RH_{\fa}^{i}\left(M\right)$  is surjective for all
$0\leq i \leq n-1$.
\end{enumerate}
\end{theorem}

\begin{prf}  (i)$\Longrightarrow$(ii) is obvious.

(ii)$\Longrightarrow$(iii) Assume that $\{x_1, \ldots, x_s\}$ satisfies the given assumptions in (ii) and set $\underline{x}:=x_{1},\ldots, x_{s}$.
Let $1\leq i\leq s$. Since $$0:_{M}\fa ^{2} \subseteq 0:_{M}x_{i}^{2}\subseteq 0:_{M}\fa \subseteq 0:_{M}\fa ^{2},$$ we have $0:_{M}\fa =0:_{M}\fa ^{2}$.
It implies that $$\RH_{\fa}^{0}\left(M\right)=0:_{M}\fa=0:_{M}\langle \underline{x}\rangle =\RH^{0}\left(\underline{x} ,M\right).$$ Thus, the natural map
$ \lambda_{M}^{0} $ is an isomorphism.

We use induction on $ n $. For $ n=1 $, we are done. Let $ n>1 $ and the result has been proved for $ n-1 $. It is enough to show that
$ \lambda_{M}^{i} $ is surjective for all $1 \leq i  \leq n-1$.

First, assume that  $\grade\left(\fa ,M\right)\geq 1$. In particular, $\Gamma_{\fa}\left(M\right)=0 $. Since $x_{1}^{2},\ldots ,x_{n}^{2}$ is
an $\fa$-weak sequence on $M$ in $ \fa ^{2} $,  Theorem \ref{2.6} implies that $\fa \RH_{\fa}^{i}\left(M\right)=0$ for all $ 0\leq i
\leq n-1 $. In particular, $ x_{1}\RH_{\fa}^{i}\left(M\right)=0 $ for all $ i<n $. Also, we have $ x_{1} \RH^{i}\left(\underline{x},M\right)=0$
for all $i \geq 0$. Since $x_{1}$ is an $\fa$-filter regular sequence on $M$, $x_{1}$ is not a zero divisor on $M$. The
exact sequence $$ 0 \longrightarrow M\stackrel{x_{1}}\longrightarrow M \longrightarrow M/x_{1}M \longrightarrow 0 $$ yields the following
commutative diagram with exact rows
\begin{displaymath}
\xymatrix{0 \ar[r] & \RH^{i-1}\left(\underline{x},M\right) \ar[r] \ar[d] & \RH^{i-1}\left(\underline{x},M/x_{1}M\right) \ar[r] \ar[d] &
\RH^{i}\left(\underline{x},M\right) \ar[r] \ar[d] & 0\\
0 \ar[r] & \RH_{\fa}^{i-1}\left(M\right) \ar[r] & \RH_{\fa}^{i-1}\left(M/x_{1}M\right) \ar[r] & \RH_{\fa}^{i}\left(M\right) \ar[r] & 0}
\end{displaymath}
for all $i<n$,  where the vertical maps are the canonical homomorphisms. Set $\overline{R}:=R/x_{1}R$, $\overline{\fa}:=\fa \overline{R}$ and
$\overline{M}:=M/x_{1}M $. We write $ \overline{x} $ for the image of $ x \in R $ in $ \overline{R} $. Since $ \overline{x}_{2},\ldots ,
\overline{x}_{s} $ is an unconditioned $ \overline{\fa} $-filter regular sequence on $ \overline{M} $ and every subset of $ \lbrace \overline{x}_{2},
\ldots, \overline{x}_{s} \rbrace $ with $ n-1 $ elements is an unconditioned $ \overline{\fa} $-weak sequence on  $ \overline{M} $, by induction
hypothesis, the natural map $ \RH^{i}\left(\overline{x}_{2},\ldots ,\overline{x}_{s}, \overline{M}\right) \longrightarrow \RH_{\overline{\fa}}^{i}\left(
\overline{M}\right) $ is surjective for all $ i\leq n-2 $.  The $R$-modules $\RH^{i}\left(\overline{x}_{2},\ldots ,\overline{x}_{s},\overline{M}\right)$
and $\RH^{i}\left(x_{2},\ldots, x_{s},\overline{M}\right)$  are isomorphic. Also, by the Independence Theorem for local cohomology, the two $R$-modules
$\RH_{\langle \overline{x}_{2},\ldots ,\overline{x}_{s} \rangle }^{i}\left(\overline{M}\right)$ and $\RH_{\langle  x_{2},\ldots, x_{s}\rangle}^{i}
\left(\overline{M}\right)$ are isomorphic.
Therefore, the corresponding natural homomorphisms $$ \RH^{i}\left(x_{2},\ldots ,x_{s},\overline{M}\right)  \longrightarrow \RH_{\langle x_{2},
\ldots, x_{s}\rangle}^{i}\left(\overline{M}\right)$$ are surjective for all $ i\leq n-2 $. By \cite[Corollary 1.7]{StV}, one has the following
commutative diagram with exact rows
\begin{displaymath}
\xymatrix{ \RH^{i}\left(\underline{x},\overline{M}\right) \ar[r] \ar[d] & \RH^{i}\left(x_{2},\ldots ,x_{s},\overline{M}\right) \ar[r] \ar[d] & 0\\
\RH_{\langle \underline{x} \rangle }^{i}\left(\overline{M}\right) \ar[r]^{\cong} & \RH_{\langle x_{2},\ldots ,x_{s} \rangle }^{i}\left(
\overline{M}\right) \ar[r] & 0}
\end{displaymath}
for all $i\geq 1$. (Note that the assumption $R$ being local is not needed in the proof \cite[Corollary 1.7]{StV}.) It follows that $\RH^{i}
\left(\underline{x},\overline{M}\right) \longrightarrow \RH_{\fa}^{i}\left(\overline{M}\right)$ is surjective for all $ i\leq n-2 $.

Now, the first commutative diagram implies that $$\RH^{i}\left(\underline{x},M\right) \longrightarrow \RH_{\fa}^{i}\left(M\right)$$ is
surjective for all $i\leq n-1$.

Now, assume that $\grade\left(\fa ,M\right)=0$.
Since $\underline{x}$ is an  unconditioned $\fa$-filter regular sequence on $M$, it follows that  $\underline{x}$ is an unconditioned $\fa$-filter
regular sequence on $ M/\Gamma_{\fa}\left(M\right)$. We show that every subset of $\lbrace x_{1},\ldots, x_{s} \rbrace$ with $n$ elements is an unconditioned
$\fa$-weak sequence on $M/ \Gamma_{\fa}\left(M\right)$. Let $\lbrace y_{1},\ldots ,y_{n} \rbrace \subseteq \lbrace x_{1},\ldots, x_{s} \rbrace$  and let
$t_{1}, \ldots ,t_{n} $ be positive integers. One has to show that $$ \left(\langle y_{1}^{t_{1}},\ldots, y_{i-1}^{t_{i-1}} \rangle M + \Gamma_{\fa}
\left(M\right) \right):_{M} y_{i}^{t_{i}}
\subseteq \left(\langle y_{1}^{t_{1}},\ldots, y_{i-1}^{t_{i-1}} \rangle M + \Gamma_{\fa}\left(M\right)\right) :_{M} \fa$$ for every $1 \leq i \leq n$.
Let $m \in \left(\langle y_{1}^{t_{1}},\ldots ,y_{i-1}^{t_{i-1}} \rangle M + \Gamma_{\fa}\left(M\right)\right) :_{M} y_{i}^{t_{i}}$. Then, there exists
$\ell \geq 1$ such that $ \fa ^{\ell} \left(y_{i}^{t_{i}}m\right) \subseteq \langle y_{1}^{t_{1}},\ldots ,y_{i-1}^{t_{i-1}} \rangle M $. So, $$ \fa ^{\ell} m
\subseteq \langle y_{1}^{t_{1}},\ldots ,y_{i-1}^{t_{i-1}} \rangle M :_{M} y_{i}^{t_{i}} \subseteq \langle y_{1}^{t_{1}},\ldots ,y_{i-1}^{t_{i-1}} \rangle
M:_{M}\fa. $$ Hence, $m\in \langle y_{1}^{t_{1}},\ldots, y_{i-1}^{t_{i-1}} \rangle M:_M\fa ^{\ell+1}$. On the other hand, we have
\begin{equation*}
\begin{array}{ll}
\langle y_{1}^{t_{1}},\ldots ,y_{i-1}^{t_{i-1}} \rangle M :_{M} \fa ^{\ell+1} & \\
 \subseteq  \langle y_{1}^{t_{1}},\ldots ,y_{i-1}^{t_{i-1}} \rangle M :_{M} y_{i}^{\ell+1} & \\
\subseteq \langle y_{1}^{t_{1}},\ldots ,y_{i-1}^{t_{i-1}} \rangle M :_{M} \fa & \\
 \subseteq \left(\langle y_{1}^{t_{1}},\ldots ,y_{i-1}^{t_{i-1}} \rangle M + \Gamma_{\fa}\left(M\right)\right) :_{M} \fa .
\end{array}
\end{equation*}
As $\grade\left(\fa, M/\Gamma_{\fa}\left(M\right)\right)>0$, the natural map $\lambda_{\frac{M}{\Gamma_{\fa}
\left(M\right)}}^{i}$  is surjective for all $i\leq n-1.$ The exact sequence $$0\rightarrow  \Gamma_{\fa}\left(M\right) \rightarrow M
\rightarrow M/\Gamma_{\fa}\left(M\right)\rightarrow 0$$ implies the following commutative diagram
\begin{displaymath}
\xymatrix{ \RH^{i}\left(\underline{x},M\right) \ar[r]^{\pi ^{i}} \ar[d]_{\lambda_{M}^{i}} & \RH^{i}\left(\underline{x},
M/\Gamma_{\fa}\left(M\right)\right) \ar[r]  \ar[d]^{\lambda_{\frac{M}{\Gamma_{\fa}\left(M\right)}}^{i}} & \RH^{i+1}\left(\underline{x},
\Gamma_{\fa}\left(M\right)\right) \ar[r]^{\mu^{i+1}} & \RH^{i+1}\left(\underline{x}, M\right)    \\
\RH_{\fa}^{i}\left(M\right) \ar[r]^{\cong}  & \RH_{\fa}^{i}\left(M/\Gamma_{\fa}\left(M\right)\right) }
\end{displaymath}
with exact top row for all $ i\geq 1 $. By Remark \ref{2.3}(A)(iii) every unconditioned $\fa$-weak sequence is $d$-sequence.  By Lemma
\ref{2.5}, $\mu^{i}$ is injective for all $i\leq n$. So, $\pi^{i}$ is surjective for $i\leq n-1$. Hence,
from the above commutative diagram, we get that $\lambda_{M}^{i}$ is surjective for all $1\leq i \leq n-1$. Since $\lambda_{M}^{0}$ is an
isomorphism, $\lambda_{M}^{i}$ is surjective for all $0 \leq i \leq n-1$.

(iii)$\Longrightarrow$(i) Assume that $\{w_1, \ldots, w_t\}$ satisfies the given assumptions in (iii). Let $\{x_1, \ldots, x_s\}$ be an
arbitrary $\text{f}$-generating set of $\fa$ with respect to $M$ with $s\geq n$  and set $\underline{x}:=x_{1},\ldots, x_{s}$. By Lemma \ref{2.1}, for each $j\geq 0$, we have the following commutative diagram in which all maps are natural
\begin{displaymath}
\xymatrix{ \RH^{j}\left(w_1, \ldots, w_t,M\right) \ar[r] \ar[d] & \RH^{j}\left(\underline{x},M\right)  \ar[d]\\
\RH_{ \langle w_1, \ldots, w_t \rangle }^{j}\left(M\right) \ar[r]^{\cong}  & \RH_{\langle \underline{x} \rangle }^{j}
\left(M\right).}
\end{displaymath}
It yields that the natural map $\lambda_{M}^{j}:\RH^{j}\left(\underline{x},M\right) \longrightarrow \RH_{\fa}^{j}\left(M\right)$
is surjective for all $0\leq j \leq n-1$.

Let $\lbrace z_{1},\ldots ,z_{n}\rbrace  \subseteq \lbrace x_{1}, \ldots, x_{s} \rbrace$ and $t_{1},\ldots, t_{n}$
be positive integers. Set $y_{i}:=z_{i}^{t_{i}}$ for all $1 \leq i \leq n $. Let $1 \leq i \leq n $ and set $X:= \lbrace x_{1}, \ldots, x_{s}
\rbrace \setminus \lbrace y_{1},\ldots, y_{i-1} \rbrace $. We choose $\alpha \in X $. As above for each $j\geq 0$, we have the following
commutative diagram
\begin{displaymath}
\xymatrix{ \RH^{j}\left(\underline{x},M\right) \ar[r] \ar[d] & \RH^{j}\left(y_{1},\ldots ,y_{i-1},\alpha ,M\right)  \ar[d]\\
\RH_{ \langle \underline{x} \rangle }^{j}\left(M\right) \ar[r]  & \RH_{ \langle y_{1},\ldots ,y_{i-1}, \alpha \rangle }^{j}
\left(M\right).}
\end{displaymath}
By Remark \ref{2.3}(B)(i),  $$\RH_{\langle \underline{x} \rangle }^{j}\left(M\right) \cong \RH_{\fa}^{j}\left(M\right)\cong
\RH_{ \langle y_{1},\ldots ,y_{i-1}, \alpha \rangle }^{j}\left(M\right)$$ for all $0 \leq  j\leq i-1$. Thus from the above diagram, we see that
the natural map $$\RH^{j}\left(y_{1},\ldots, y_{i-1},\alpha ,M\right) \longrightarrow \RH^{j}_{ \langle y_{1},\ldots ,y_{i-1}, \alpha \rangle }
\left(M\right)$$ is surjective for all $0 \leq j\leq i-1 $. Hence by Lemma \ref{2.4},  $y_{1},\ldots ,y_{i-1},\alpha $ is an $u.s.d$-sequence on $M$. Set $\overline{M}:=M/ \langle y_{1},\ldots, y_{i-1} \rangle M$. As $$\langle y_{1},\ldots, y_{i-1} \rangle M:_{M}\alpha=\langle y_{1},
\ldots, y_{i-1} \rangle M:_{M}\alpha^{2},$$ we get $0:_{\overline{M}}\alpha =0:_{\overline{M}}\alpha^{2}$. Thus, as $\alpha$ is an $\fa$-filter regular sequence on $\overline{M}$, one has $$\Gamma_{\fa}\left(\overline{M}\right)=\Gamma_{\langle \alpha \rangle}\left(\overline{M}\right)=0:_{\overline{M}}\alpha.$$ Now, we have
\begin{align*}
0:_{\overline{M}}y_{i} &\subseteq \Gamma_{ \langle y_{i} \rangle }\left(\overline{M}\right)\\
&=\Gamma_{\fa}\left(\overline{M}\right)\\
&=\bigcap_{\alpha \in X}\left(0:_{\overline{M}}\alpha\right)\\
&=0:_{\overline{M}}\fa,
\end{align*}
and so $\langle y_{1},\ldots ,y_{i-1} \rangle M:_{M}y_{i}\subseteq \langle y_{1},\ldots ,y_{i-1} \rangle M:_{M}\fa$.
\end{prf}

\section{Relative generalized Cohen-Macaulay modules}

The class of generalized Cohen-Macaulay modules contains the class of Cohen-Macaulay modules. Indeed generalized Cohen-Macaulay
modules enjoy many interesting properties similar to the ones of Cohen-Macaulay modules. As a generalization of the notion of
Cohen-Macaulay modules, in \cite{RZ} and \cite{DGTZ}, relative Cohen-Macaulay modules were studied. As a continuation of the
work, we establish a theory of relative generalized Cohen-Macaulay modules. To this end, first we specify some notation and facts.

\begin{remark} \label{3.1} Let $\fa$ be an ideal of $R$ and $M$ a finitely generated $R$-module.
\begin{enumerate}
\item[(i)] Following \cite[Definitions 9.2.2]{BS}, we set $\lambda_{\fa}\left(M\right):=\inf \lbrace \depth M_{\fp}+
\Ht\left(\frac{\fa +\fp}{\fp}\right)|  \  \ \fp\in \Spec R \setminus \V\left(\fa\right) \rbrace.$ By \cite[Theorem 9.3.7]{BS}, we
have $\text{f}_{\fa}\left(M\right) \leq \lambda _{\fa}\left(M\right)$.
\item[(ii)] Faltings' Local-global Principle Theorem \cite[Satz 1]{F} implies that $$\text{f}_{\fa}\left(M\right)=\inf \lbrace i\in
\mathbb{N}|\ \fa^{\ell}\RH_{\fa}^{i}\left(M\right)\neq 0 \ \text{for\ all}\ {\ell}\in \mathbb{N} \rbrace.$$
\item[(iii)] If $\fa M \neq M$, then $\Ht_{M}\fa \leq \cd\left(\fa ,M\right)$, otherwise one has  $\text{ht}_{M} \fa=+\infty$
and $\cd\left(\fa,M\right)=-\infty $. (Recall that, by convention, $\inf \emptyset=+\infty$ and $\sup \emptyset=-\infty $.)
\item[(iv)] If $c:=\cd(\fa,M)>0$, then by \cite[Corollary 3.3(i)]{DV}, the $R$-module $\RH_{\fa}^{c}\left(M\right)$ is not finitely
generated. So in this case, one has $\text{f}_{\fa}\left(M\right)\leq \cd(\fa,M)$. From this, we can also deduce that if
$\Ht_{M} \fa>0$, then $\text{f}_{\fa}\left(M\right)\leq \Ht_{M} \fa$.
\end{enumerate}
\end{remark}

\begin{definition}\label{3.11} Let $\fa$ be an ideal of $R$. A finitely generated $R$-module $M$ is said to be
$\fa$-\it{relative generalized Cohen-Macaulay} if $\cd\left(\fa,M\right) \leq 0$; or $\cd\left(\fa,M\right)=
\text{f}_{\fa}\left(M\right)$.
\end{definition}

The next lemma is needed in the proof of Theorem \ref{3.3}. Theorem \ref{3.3} is one of the main results of the paper which
provides some properties of relative system of parameters whenever the module under consideration is relative generalized
Cohen-Macaulay.

\begin{lemma}\label{3.2} Let $\fa$ be an ideal of $R$ which is contained in the Jacobson radical of $R$ and $M$ an $\fa$-relative generalized Cohen-Macaulay
$R$-module with $c:=\cd\left(\fa,M\right)>0$. Then $\text{f}_{\fa}\left(M\right)=\lambda_{\fa}\left(M\right)$ and $\cd\left(\fa,R/\fp\right)=c$ for every
$\fp \in \Ass_{R} M \setminus \V\left(\fa\right)$.
\end{lemma}

\begin{prf}  First of all note that  the assumption $ c>0 $, implies that $ \Ass_{R} M \nsubseteq \V\left(\fa\right) $. Let $ \fp \in \Ass_{R} M \setminus
\V\left(\fa\right)$.  Then $ \text{depth}_{R_{\fp}}M_{\fp}=0 $. Using the assumption,  \cite[Theorem 9.3.7]{BS} and \cite[Theorem 2.2]{DNT} yields that
\begin{equation*}
\begin{array}{ll}
\cd\left(\fa ,M\right)=\text{f}_{\fa}\left(M\right) &\\
\leq \lambda_{\fa}\left(M\right) & \\
\leq \text{depth}_{R_{\fp}}M_{\fp} +\text{ht}\left(\frac{\fa + \fp}
{\fp}\right) &\\
=\text{ht}\left(\frac{\fa +\fp}{\fp}\right) &\\
\leq  \cd\left(\frac{\fa +\fp}{\fp}, R/ \fp\right) & \\
=\cd\left(\fa ,R/ \fp\right) &\\
\leq \cd\left(\fa ,M\right).

\end{array}
\end{equation*}
\end{prf}

\begin{theorem}\label{3.3} Let $\fa$ be an ideal of $R$ which is contained in the Jacobson radical of $R$ and $M$ an $\fa$-relative generalized Cohen-Macaulay
$R$-module such that $c:=\cd\left(\fa ,M\right)=\ara\left(\fa ,M\right)$. Assume that $\underline{x}:= x_{1}, \ldots, x_{c}$ is an $\fa$-s.o.p of M. Then
\begin{enumerate}
\item[(i)] $x_{1},\ldots ,x_{c}$ is an unconditioned $\fa$-filter regular sequence on $M$.
\item[(ii)] $M/ \langle x_{1},\ldots, x_{i} \rangle M$ is an $\fa$-relative generalized Cohen-Macaulay $R$-module and $$\cd\left(\fa,M/\langle x_{1},
\ldots, x_{i} \rangle M\right)=c-i$$ for all $i=0,\ldots, c$.
\end{enumerate}
\end{theorem}

\begin{prf} Clearly if $z_{1}, z_2, \ldots, z_{c}\in \fa$ is an $\fa$-s.o.p of $M$, then for all $t_{1},\ldots,t_{c}\in \mathbb{N}$, every permutation of
$z_{1}^{t_{1}},\ldots, z_{c}^{t_{c}}$ is also an $\fa$-s.o.p of $M$. So to prove (i), it is enough to show that $x_{1},\ldots, x_{c}$ is an $\fa$-filter regular
sequence on $M$. In case $c=0$, there is nothing to prove. Hence, we may assume that $c>0$. It is enough to show that the claim holds for $i=1$. Set $x:=x_{1}.$

(i) To contrary, assume that there is $ \fp \in \Ass_{R} M \setminus \V\left(\fa\right) $ such that $x \in \fp$. Note that $$\Rad \left( \Ann_{R}
\left(M/\fp M\right)\right)=\Rad\left(\fp+\Ann_{R} M\right)=\Rad \fp =\fp. $$ Thus $ \Supp_{R}\left(M/\fp M\right)=\Supp_{R}\left(R/ \fp\right) $, and so
by \cite[Theorem 2.2]{DNT}, $\cd\left(\fa,R/ \fp\right)=\cd\left(\fa,M/ \fp M\right)$. So by Lemma \ref{3.2}, $\cd\left(\fa ,M/ \fp M\right)=c$.
As $x\in \fp$, one has a natural $R$-epimorphism $M/xM \longrightarrow M/\fp M$, and hence $\Supp_{R}\left(M/ \fp M\right)\subseteq \Supp_{R}
\left(M/xM\right)$, by employing \cite[Theorem 2.2]{DNT} again. Therefore, $$c=\cd\left(\fa , M/ \fp M\right) \leq \cd\left(\fa , M/xM\right).$$
On the other hand, by \cite[ Lemma 2.4]{DGTZ}, $\cd\left(\fa,M/xM\right)=c-1$,  which yields our desired contradiction.

(ii) Now, we prove that $M/xM$ is an  $\fa$-relative generalized Cohen-Macaulay $R$-module. One has $ \cd\left(\fa ,M/xM\right)=c-1 $. So if $c=1$,
the claim holds. Hence, we may assume that $ c>1 $. Note that  $$\text{f}_{\fa}\left(M/xM\right) \leq \cd\left(\fa ,M/xM\right)=c-1.$$ Thus, to complete
the proof, we need to show  that $ \RH_{\fa}^{i}\left(M/xM\right) $ is finitely generated for all  $ 0 \leq i \leq c-2 $. Since $\RH_{\fa}^{0}\left(M/xM
\right)$ is finitely generated, we can assume that $ 1 \leq i \leq c-2 $. Since  $ x $ is an $\fa$-filter regular sequence on $M$,  $x$
is non-zero divisor on $ M/ \Gamma_{\fa}\left(M\right) $. Applying the functor $ \RH_{\fa}^{i}\left(-\right) $ on the exact sequence
$$ 0 \longrightarrow M/ \Gamma_{\fa}\left(M\right) \stackrel{x}\longrightarrow M/ \Gamma_{\fa}\left(M\right) \longrightarrow M/\left(xM+\Gamma_{\fa}
\left(M\right)\right) \longrightarrow 0, $$
yields the following exact sequence
$$ \RH_{\fa}^{i}\left(M/ \Gamma_{\fa}\left(M\right)\right) \longrightarrow \RH_{\fa}^{i}\left(M/\left(xM+ \Gamma_{\fa}\left(M\right)\right)\right)
\longrightarrow \RH_{\fa}^{i+1}\left(M/
\Gamma_{\fa}\left(M\right)\right). $$
Since $ \RH_{\fa}^{j}\left(M/ \Gamma_{\fa}\left(M\right)\right)\cong \RH_{\fa}^{j}\left(M\right) $ is finitely generated for $ j=i,i+1 $, the above
exact sequence implies
that $ \RH_{\fa}^{i}\left(M/\left(xM+\Gamma_{\fa}\left(M\right)\right)\right) $ is finitely generated. From the exact sequence
$$ 0 \longrightarrow \left(xM+\Gamma_{\fa}\left(M\right)\right)/xM \longrightarrow M/xM \longrightarrow M/\left(xM+\Gamma_{\fa}\left(M\right)\right)
\longrightarrow 0,$$ we obtain the following
exact sequence
$$ \RH_{\fa}^{i}\left( \left(xM+\Gamma_{\fa}\left(M\right)\right)/xM\right) \longrightarrow \RH_{\fa}^{i}\left(M/xM\right) \longrightarrow \RH_{\fa}^{i}\left(M/\left(xM+\Gamma_{\fa}\left(M\right)\right)\right) $$
As $\left(xM+\Gamma_{\fa}\left(M\right)\right)/xM $ is $\fa$-torsion,  it follows that $\RH_{\fa}^{i}\left(\left(xM+\Gamma_{\fa}
\left(M\right)\right)/xM\right)=0 $. So, the above exact sequence implies
that $ \RH_{\fa}^{i}\left(M/xM\right) $ is finitely generated.
\end{prf}

Part (ii) of the following example provides an ideal $\fa$ and an $\fa$-relative generalized Cohen-Macaulay $R$-module $M$ such that $\fa$ contains
no $\fa$-s.o.p of $M$.

\begin{example}\label{3.4}
\begin{enumerate}
\item[(i)] Let $\fa$ be an ideal of $R$. Any finitely generated $R$-module with $\cd\left(\fa,M\right)=1$ is an $\fa$-relative
generalized Cohen-Macaulay $R$-module.
\item[(ii)]  Let $K$ be a field and $S:=K[[X,Y,Z,W]]$. Consider the elements $f:=XW-YZ$, $g:=Y^{3}-X^{2}Z$ and $h:=Z^{3}-Y^2W$
of $S$ and set $R:=S/\langle f \rangle$ and $\fa:=\langle f,g,h \rangle/\langle f \rangle$. Then $R$ is a Noetherian local ring
of dimension $3$, $\cd\left(\fa,R\right)=1$ and $\ara\left(\fa,R\right)\geq 2$; see \cite[Remark 2.1(ii)]{HS}.  As
$\ara\left(\fa,R\right)\neq \cd\left(\fa,R\right)$, it follows that $\fa$ contains no $\fa$-s.o.p of $R$. But, by (i) $R$ is an
$\fa$-relative generalized Cohen-Macaulay.
\end{enumerate}
\end{example}

The following result provides a characterization of $\fa$-relative generalized Cohen-Macaulay modules by using $\fa$-weak sequences.

\begin{proposition}\label{3.5} Let $\fa$  be an ideal of $R$ which is contained in the Jacobson radical of $R$ and $M$ a finitely generated $R$-module.
Suppose that $\cd \left(\fa,M\right)= \ara \left(\fa,M\right) $. Then the following conditions are equivalent:
\begin{enumerate}
\item[(i)] $M$ is $\fa$-relative generalized Cohen-Macaulay.
\item[(ii)] There is an integer $\ell$ such that every $\fa$-s.o.p of $M$ is an $\fa ^{\ell}$-weak sequence on $M$.
\end{enumerate}
\end{proposition}

\begin{prf} (i)$\Longrightarrow$(ii) By Theorem \ref{3.3}, every $\fa$-s.o.p of $M$ is an $\fa$-filter regular sequence on $M$. So, the assertion follows by
\cite[Theorem (i)$\Rightarrow$(iii)]{KS}.

(ii)$\Longrightarrow$(i) It follows by \cite[Theorem (iv)$\Rightarrow$(i)]{KS}.
\end{prf}

\begin{proposition}\label{3.6}  Let $\fa$ be an ideal of $R$   and $M$ be an $\fa$-relative generalized Cohen-Macaulay
$R$-module with $ \cd\left(\fa ,M\right)>0 $. Then
  $M_{\fp}$ is Cohen-Macaulay for every $\fp\in \Supp_{R} M $ with $ \fp \subsetneqq \fa $.
\end{proposition}

\begin{prf} Let $\fp \in \Supp_{R} M$ such that  $\fp \subsetneqq \fa$.  Then $\fp \in \Supp_{R} M \setminus \V\left(\fa\right)$. One has
\begin{equation*}
\begin{array}{ll}
\cd\left(\fa,M\right)=\text{f}_{\fa}\left(M\right) &\\
\leq \lambda_{\fa}\left(M\right)  & \\
\leq \text{depth}_{R_{\fp}}M_{\fp} + \text{ht}\left(\left(\fa +\fp\right)/\fp\right)  & \\
=\text{depth}_{R_{\fp}}M_{\fp} + \text{ht}\left(\fa / \fp\right) &\\
\leq \text{dim}_{R_{\fp}}M_{\fp}+\text{ht}\left(\fa / \fp\right)&\\
\leq \text{ht}_{M}\fa &\\
\leq \cd\left(\fa,M\right).

\end{array}
\end{equation*}
Thus $\depth_{R_{\fp}}M_{\fp}=\dim_{R_{\fp}}M_{\fp}$, as required.
\end{prf}

\section{Relative Buchsbaumness}

We start this section by introducing the notions of $\fa$-relative surjective Buchsbaum, $\fa$-relative Buchsbaum
and $\fa$-relative quasi Buchsbaum $R$-modules.

\begin{definition}\label{4.11} Let $\fa$ be an ideal of $R$  and $M$ a finitely generated $R$-module.
\begin{enumerate}
\item[(i)] We say that $M$ is $\fa$-\it{relative surjective Buchsbaum} if the natural map $$\varphi_{M}^{i}:\Ext_{R}^{i}
\left(R/\fa , M\right) \longrightarrow \RH_{\fa}^{i}\left(M\right)$$ is surjective for all $i<\cd\left(\fa,M\right)$.
\item[(ii)]  We say that $M$ is $\fa$-\it{relative Buchsbaum} if $\cd\left(\fa,M\right)=\ara\left(\fa ,M\right)$ and every
$\fa$-s.o.p of $M$ is an $\fa$-weak sequence on $M$.
\item[(iii)]  We say that $M$ is $\fa$-\it{relative quasi Buchsbaum} if $\cd\left(\fa ,M\right)=\ara\left(\fa ,M\right)$
and every $\fa$-s.o.p of $M$ contained in $\fa^2$ is an $\fa$-weak sequence on $M$.
\end{enumerate}
\end{definition}

Next, we establish the following characterization of relative quasi Buchsbaum modules.

\begin{theorem}\label{4.1} Let $\fa$ be an ideal of $R$ and $M$ a finitely generated $R$-module with $\cd\left(\fa ,M\right)=\ara\left(\fa ,M\right)$.
Consider the following conditions:
\begin{enumerate}
\item[(i)] $M$ is $\fa$-relative quasi Buchsbaum.
\item[(ii)] There is an $\fa$-s.o.p of $M$ contained in $\fa ^{2}$ which is an $\fa$-weak sequence on $M$.
\item[(iii)]  $\fa\RH_{\fa}^{i}\left(M\right)=0$ for all $0 \leq i<\cd\left(\fa,M\right)$.
\end{enumerate}
Then $\left(i\right)$ implies $\left(ii\right)$ and $\left(ii\right)$ implies $\left(iii\right)$. Furthermore if $\fa$ is
contained in the Jacobson radical of $R$, then $\left(iii\right)$ implies $\left(i\right)$.
\end{theorem}

\begin{prf} (i)$\Longrightarrow$(ii) holds by the definition.

(ii)$\Longrightarrow$(iii) follows by Theorem \ref{2.6}.

(iii)$\Longrightarrow$(i) Assumption (iii) implies that $M$ is  $\fa $-relative generalized Cohen-Macaulay. Let $ x_{1},\ldots, x_{c}\in \fa ^{2}$
be an $\fa$-s.o.p of $M$. By Theorem \ref{3.3}, the sequence $ x_{1},\ldots ,x_{c} $ is an  $\fa$-filter regular sequence on $M$. So, by Theorem
\ref{2.6}, it turns out that $x_{1},\ldots, x_{c}$ is an $\fa$-weak sequence on $M$.
\end{prf}

Goto in \cite[Corollary 2.8]{G} proved that a finitely generated $R$-module $M$ is Buchsbaum if and only if every system of parameters of $M$ forms
a $d$-sequence on $M$. Accordingly, the authors tried to show that the conditions (iii) and (iv) in the following result are equivalent, but so
far without success.

\begin{theorem}\label{4.2} Let $\fa$ be an ideal of $R$ and $M$ a finitely generated $R$-module.  Assume that $ c:=\cd\left(\fa ,M\right)=\ara\left(\fa
,M\right)>0  $.  Consider the following conditions:
\begin{enumerate}
\item[(i)] For every  generating set $\{b_{1},\ldots, b_{t}\}$ of $\fa$, the natural map $\lambda_{M}^{i}:\RH^{i}\left(b_{1},\ldots ,b_{t},M\right)
\longrightarrow \RH_{\fa}^{i}\left(M\right)$ is surjective for all $0\leq i \leq c-1$.
\item[(ii)] There exists an $f$-generating set $\{a_{1},\ldots, a_{s}\}$ of $\fa$ with respect to $M$ such that the natural map $\lambda_{M}^{i}:\RH^{i}
\left(a_{1},\ldots, a_{s},M\right)\longrightarrow \RH_{\fa}^{i}\left(M\right)$ is surjective for all $0\leq i \leq c-1$.
\item[(iii)] $M$ is $\fa$-relative Buchsbaum.
\item[(iv)]  Every $\fa$-s.o.p of $M$  is an $u.s.d$-sequence on $M$.
\end{enumerate}
Then $\left(i\right)$ and  $\left(ii\right)$ are equivalent and  $\left(iii\right)$ implies  $\left(iv\right)$. Furthermore if $\fa$ is contained
in the Jacobson radical of $R$, then $\left(i\right)$ implies $ \left(iii\right) $.
\end{theorem}

\begin{prf} (i)$\Longrightarrow$(ii) is clear by Remark \ref{2.3}(B)(ii).

(ii)$\Longrightarrow$(i) Let $\{b_{1},\ldots, b_{t}\}$ be a generating set of $\fa$. By Lemma \ref{2.1}, one has the following commutative diagram
\begin{displaymath}
\xymatrix{ \RH^{i}\left(a_{1},\ldots ,a_{s},M\right) \ar[r] \ar[d] & \RH^{i}\left(b_{1},\ldots ,b_{t},M\right)  \ar[d]   \\
  \RH_{\fa}^{i}\left(M\right) \ar[r]^{\cong}  & \RH_{\fa}^{i}\left(M\right) }
\end{displaymath}
for all $ 0 \leq i \leq c-1 $. It implies that $ \RH^{i}\left(b_{1},\ldots ,b_{t},M\right) \longrightarrow \RH_{\fa}^{i}\left(M\right) $ is surjective for all
$ 0 \leq i \leq c-1 $.

(iii)$\Longrightarrow$(iv) follows by Remark \ref{2.3}(A)(iii).

(i)$\Longrightarrow$ (iii) Let $ x_{1},\ldots ,x_{c} \in \fa $ be an $\fa$-s.o.p of $M$. The assumption implies that  $ \fa \RH_{\fa}^{i}\left(M\right)=0 $
for $ 0 \leq i \leq c-1 $. Hence by Theorem \ref{3.3}, $ x_{1},\ldots ,x_{c} $ is an unconditioned $\fa$-filter regular sequence on $M$. We can find $ z_{1},
\ldots, z_{l}$ in $\fa$ such that $\fa=\langle x_{1},\ldots ,x_{c},z_{1},\ldots, z_{l} \rangle $. As in the proof of \cite[Proposition 1.2]{TZ}, there
are $w_{1},\ldots, w_{l}$ in $\fa$ such that $\fa=\langle x_{1},\ldots ,x_{c},w_{1},\ldots, w_{l}\rangle$ and $x_{1},\ldots, x_{c}, w_{1},\ldots,
w_{l}$ is an unconditioned $\fa$-filter regular sequence on $M$. By the hypothesis, the natural map $\RH^{i}\left(x_{1},\ldots ,x_{c},w_{1},\ldots ,w_{l}\right)
\longrightarrow \RH_{\fa}^{i}\left(M\right)$ is surjective for all $0 \leq i \leq c-1$.   The claim follows from Theorem \ref{2.7}.
\end{prf}

Next, we provide the comparison between the classes of modules that are introduced in Definitions \ref{3.11} and \ref{4.11}.

\begin{corollary}\label{4.3}
\begin{enumerate}
\item[(i)] Any $\fa$-relative Cohen-Macaulay $R$-module is $\fa$-relative surjective Buchsbaum.
\item[(ii)] If the ideal $\fa$ is contained in the Jacobson radical of $R$ and $M$ is an $\fa$-relative surjective Buchsbaum $R$-module with
$\cd\left(\fa,M\right)=\ara\left(\fa,M\right)$, then $M$ is $\fa$-relative Buchsbaum.
\item[(iii)] Any $\fa$-relative Buchsbaum $R$-module is $\fa$-relative quasi Buchsbaum.
\item[(iv)] Any $\fa$-relative quasi Buchsbaum $R$-module is $\fa$-relative generalized Cohen-Macaulay.
\end{enumerate}
\end{corollary}

\begin{prf} (i) and (iii) are obvious by the definitions.

(ii) Let $\fa=\langle x_{1},\ldots, x_{s} \rangle$. By \cite[Lemma 1.5]{StV},  we have the following commutative diagram

\begin{displaymath}
 \xymatrix{\Ext_{R}^{i}\left(R/\fa ,M\right) \ar[r]^{\varphi_{M}^{i}}  \ar[d]_{\psi_{M}^{i}}  & \RH_{\fa}^{i}\left(M\right) \\
\RH^{i}\left(x_{1},\ldots ,x_{s} ,M\right) \ar[ur]_{\lambda_{M}^{i}}}
\end{displaymath}
Since $\varphi_{M}^{i}$ is surjective for all $i<c:=\cd(\fa,M)$ and $\lambda_{M}^{i}\psi_{M}^{i}=\varphi_{M}^{i}$, we deduce
that $\lambda_{M}^{i}$ is surjective for all $i<c$. Thus the claim follows by Theorem \ref{4.2}.

(iv) Theorem \ref{4.1} implies that $\fa\RH_{\fa}^{i}\left(M\right)=0$ for all $i<\cd\left(\fa,M\right)$. Then by Faltings'
Local-global Principle Theorem \cite[Satz 1]{F}, $\RH_{\fa}^{i}\left(M\right)$ is finitely generated for all $i<\cd\left(\fa,
M\right)$. Thus, either $\cd\left(\fa,M\right)\leq 0$ or $\cd\left(\fa,M\right)=\text{f}_{\fa}\left(M\right)$.
\end{prf}

Below, we record another corollary of Theorem \ref{4.2}.

\begin{corollary}\label{4.4} Let $(R,\fm)$ be a local ring and $M$ a non-zero finitely generated $R$-module. Let $\fa$ be an
ideal of $R$ which is generated by a part of a system of parameters of $M$. Then $\cd(\fa,M)=\ara(\fa,M)$. Moreover, if $M$
is Buchsbaum, then it is also $\fa$-relative Buchsbaum.
\end{corollary}

\begin{prf} Let $d:=\dim M$. Suppose that $\fa:=\langle x_{1}, x_2, \ldots, x_n\rangle$, where $x_{1}, \ldots, x_{n}, \ldots, x_d$
is a system of parameters of $M$. Then $x_{1}, x_2, \ldots, x_n$ is a system of parameters for the $R$-module $N:=M/\langle x_{n+1},
\ldots, x_d\rangle M$, and so $\RH^{n}_{\fa}(N)\cong \RH^{n}_{\fm}(N)\neq 0$. Now, the short exact sequence $$0\lo \langle x_{n+1},
\ldots, x_d\rangle M\lo M\lo N\lo 0,$$ yields an $R$-epimorphism $\varphi:\RH^{n}_{\fa}(M)\lo \RH^{n}_{\fa}(N)$. Therefore,
$\RH^{n}_{\fa}(M)\neq 0$, and so $\cd(\fa,M)=n=\ara(\fa,M)$.

Next, assume that $M$ is Buchsbaum. By the definition of Buchsbaum modules, it turns out that $x_{1}, \ldots, x_{n}$ is an unconditioned
$\fm$-weak sequence on $M$. Hence, by  Remark \ref{2.3}(A)(iii), it follows  that $x_{1},\ldots, x_{n}$ is an $u.s.d$-sequence on $M$.
So, by \cite[Theorem 2.6]{G}, the canonical map $$\RH^{i}( x_{1},\ldots, x_{n},M)\longrightarrow \RH^{i}_{\langle  x_{1},\ldots,
x_{n}\rangle}(M)$$ is surjective for all $0\leq i\leq n-1$. Also, Remark \ref{2.3}(A)(ii) yields that $x_{1}, \ldots,
x_{n}$ is an unconditioned $\fa$-filter regular sequence on $M$. Thus, by Theorem \ref{4.2}, $M$ is $\fa$-relative Buchsbaum.
\end{prf}

The following example shows that the condition $\fa$ is contained in the Jacobson radical of $R$ is necessary in Corollary \ref{4.3}(ii).

\begin{example}\label{4.5} Let $k$ be a field and $u$, $v$, $w$ be independent variables. Consider the polynomial ring $R:=k[u,v,w]$. Set
$a_{1}:= v\left(1-u\right)$, $a_{2}:=w\left(1-u\right)$, $a_{3}:=u$ and $\fa:=\langle a_{1},a_{2},a_{3}\rangle$. Note that $\fa$ is not
contained in the Jacobson radical of $R$ and $\fa=\langle u,v,w \rangle $. We have $\RH_{\fa}^{i}\left(R\right)=0$ for every $0 \leq i <3$
and $\RH_{\fa}^{3}\left(R\right)\neq 0$. So, $\cd\left(\fa,R\right)=\ara\left(\fa, R\right)=3$. Hence, $R$ is $\fa$-relative Cohen-Macaulay.
If $R$ is $\fa$-relative Buchsbaum, then by Theorem \ref{4.2} the sequence $a_{1}, a_{2}, a_{3}$ is a $d$-sequence. On the other hand, $v
\in \langle a_{1} \rangle :_{R}a_{2}a_{3}$, but $v \notin \langle a_{1}\rangle :_{R} a_{3}$. So, $a_{1},a_{2}, a_{3}$ is not a $d$-sequence
on $R$.
\end{example}

\begin{proposition}\label{4.6} Let $\fa$ be an ideal of $R$ which is contained in the Jacobson radical of $R$ and $M$ a finitely generated
$R$-module with $\cd\left(\fa,M\right)=\ara\left(\fa,M\right).$ Assume that $r:=\grade\left(\fa,M\right)<c:=\cd\left(\fa,M\right)$ and
$\RH_{\fa}^{i}\left(M\right)=0$ for all $i\neq r,c$. Then the following are equivalent:
\begin{enumerate}
\item[(i)] $M$ is $\fa$-relative surjective Buchsbaum.
\item[(ii)] $M$ is $\fa$-relative Buchsbaum.
\item[(iii)] $M$ is $\fa$-relative quasi Buchsbaum.
\end{enumerate}
\end{proposition}

\begin{prf} (i)$\Longrightarrow$(ii) and (ii)$\Longrightarrow$(iii) follow by Corollary \ref{4.3}.

(iii)$\Longrightarrow$(i) Since $\grade\left(\fa,M\right)=r$, \cite[Lemma 2.5(b)]{AS} implies the following natural $R$-isomorphism
$$\Ext_{R}^{r}\left(R/\fa, M\right)\cong \Hom_R(R/\fa,\RH_{\fa} ^{r}\left(M\right)).$$ As $M$ is $\fa$-relative quasi Buchsbaum, Theorem
\ref{4.1} implies that $\fa \RH_{\fa}^{r}(M)=0$, and so one deduces the natural $R$-isomorphism $\Ext_{R}^{r}\left(R/\fa, M\right)\cong
\RH_{\fa}^{r}\left(M\right)$. Thus, $M$ is $\fa$-relative surjective Buchsbaum.
\end{prf}

The following lemma is needed in the proof of Theorem \ref{4.8}.

\begin{lemma}\label{4.7}  Let $\fa$ be an ideal of $R$ and $M$ a finitely generated $R$-module. Let $n$ be a natural number and
$x_{1},\ldots, x_{n} \in \fa$. If $\fa \Gamma_{\langle x_{k+1} \rangle } \left(M/ \langle x_{1},\ldots, x_{k} \rangle M\right)=0$
for all $0 \leq k \leq n-1$,  then $x_{1},\ldots, x_{n}$ is both  $\fa$-weak sequence on $M$ and $d$-sequence on $M$.
\end{lemma}

\begin{prf} Let $ 0 \leq k \leq n-1 $ and set $ N:=\frac{M}{\langle x_{1},\ldots ,x_{k} \rangle M} $. We have
$$
\begin{array}{ll}
0:_{N}x_{k+1}&\subseteq 0:_{N}x_{k+1}^{2} \\
&\subseteq \Gamma_{\langle x_{k+1} \rangle}\left(N\right) \\
&=0:_{N} \fa \\
&\subseteq 0:_{N} \langle x_{1},\ldots, x_{n}\rangle \\
&\subseteq 0:_{N}x_{k+1}.
\end{array}
$$
Thus $0:_{N}x_{k+1}=0:_{N}\fa$, and so $x_{1},\ldots, x_{n}$ is an $\fa$-weak sequence on $M$.

Since $ 0:_{N}x_{k+1}^{2}=0:_{N} \langle x_{1},\ldots, x_{n} \rangle,$ by \cite[Theorem 1.1(v)]{T1}, one deduces that $x_{1},\ldots, x_{n}$
is a  $d$-sequence on  $M$.
\end{prf}

The notion of relative system of parameters plays an important role in the proof of the next result which establishes a characterization of
relative Buchsbaum modules in terms of relative quasi Buchsbaum modules.

\begin{theorem}\label{4.8}
Let $\fa$ be an ideal of $R$ which is contained in the Jacobson radical of $R$ and let $M$ be a finitely generated $R$-module with
$c:=\cd\left(\fa ,M\right)=\ara\left(\fa,M\right)$. Then the following are equivalent:
\begin{enumerate}
\item[(i)] For every $\fa$-s.o.p $x_{1},\ldots, x_{c}$ of $M$, $\fa \Gamma_{\langle x_{k+1} \rangle } \left(M/ \langle x_{1},\ldots, x_{k}
\rangle M\right)=0 $ for all $0 \leq k \leq c-1$.
\item[(ii)] $M$ is $\fa$-relative Buchsbaum.
\item[(iii)] For every $\fa$-s.o.p $ x_{1},\ldots ,x_{c} \in \fa $ and every $ 0\leq k \leq c-1 $, the $R$-module $M/ \langle x_{1},\ldots, x_{k}
\rangle M$ is $\fa$-relative quasi Buchsbaum.
\end{enumerate}
\end{theorem}

\begin{prf} (i)$\Longrightarrow$(ii) By Lemma \ref{4.7}, every $\fa$-s.o.p of $M$ is an $\fa$-weak sequence on $M$. So, $M$ is $\fa$-relative
Buchsbaum.

(ii)$\Longrightarrow$(iii) Let $M$ be an $\fa$-relative Buchsbaum $R$-module and $ x_{1},\ldots ,x_{c} $ be an  $\fa$-s.o.p of $M$. Let $0
\leq k \leq c-1$. Note that $\underline{ \acute{x}}:=x_{1},\ldots, x_{k}, x_{k+1}^{2}, \ldots, x_{c}^{2} $ is also an $\fa$-s.o.p of $M$. Hence,
$\underline{\acute{x}}$ is an  $\fa$-weak sequence on $M$. So, $x_{k+1}^{2},\ldots, x_{c}^{2}$ is an $\fa$-weak sequence on $M/\langle x_{1},
\ldots, x_{k} \rangle M$. Thus, by Theorem \ref{2.6}, $\fa \RH_{\fa}^{i}\left(M/\langle x_{1},\ldots ,x_{k} \rangle M\right)=0$ for
all $0 \leq i < c-k$. Also,  by \cite[Lemma 2.4]{DGTZ}, it turns out that $$\cd(\fa,M/\langle x_{1},\ldots, x_{k} \rangle M)=\ara(\fa,M/\langle
x_{1},\ldots, x_{k}\rangle M)=c-k.$$ Therefore, by Theorem \ref{4.1}, $M/\langle x_{1},\ldots ,x_{k} \rangle M$ is $\fa$-relative quasi-Buchsbaum.

(iii)$\Longrightarrow$(i) Let $ x_{1},\ldots, x_{c} $ be an  $\fa$-s.o.p of $M$ and let $ 0 \leq k \leq c-1 $. By the assumption and Corollary
\ref{4.3},  $M/\langle x_{1},\ldots, x_{k} \rangle M$ is $\fa$-relative generalized Cohen-Macaulay. Since by \cite[Lemma 2.4]{DGTZ},
$x_{k+1},\ldots, x_{c}$ is an $\fa$-s.o.p of $M/\langle x_{1},\ldots ,x_{k} \rangle M$,  Theorem \ref{3.3}(i) implies that  $x_{k+1},
\ldots, x_{c}$ is an  $\fa$-filter regular sequence on $M/\langle x_{1},\ldots, x_{k} \rangle M$. Therefore, by Remark \ref{2.3}(B)(i)
$$\Gamma_{\fa}\left(M/\langle x_{1},\ldots, x_{k} \rangle M\right)=\Gamma_{\langle x_{k+1} \rangle} \left(M/\langle x_{1},\ldots,
x_{k}\rangle M\right), $$ and so Theorem \ref{4.1}  yields that $\fa \Gamma_{\langle x_{k+1} \rangle} \left( M/\langle x_{1},\ldots,
x_{k} \rangle M\right)=0$.
\end{prf}

\begin{acknowledgement} The authors thank the anonymous referee, whose valuable comments greatly improved the exposition
of the paper.
\end{acknowledgement}

%%%%%%%%%%%%%%%%%%%%%%%%%%%%%%%%%%%%%%%%%%%%%%%%%%%

\end{document}